\documentclass{amsart}

\usepackage{graphicx}

\newtheorem{theorem}{Theorem}[section]

\theoremstyle{definition}

\theoremstyle{remark}
\newtheorem{remark}[theorem]{Remark}

\numberwithin{equation}{section}

\usepackage{hyperref}

\begin{document}

\title{On surface completion and image inpainting by biharmonic functions: Numerical aspects}

\author{S. B. Damelin}

\author{N. S. Hoang}

\subjclass[2000]{}
\date{}

\begin{abstract}
Numerical experiments with smooth surface extension and image inpainting using harmonic and biharmonic functions are carried out. 
The boundary data used for constructing biharmonic functions are the values of the Laplacian and normal derivatives of the functions on the boundary. 
Finite difference schemes for solving these harmonic functions are discussed in detail. 
\end{abstract}


\maketitle

\section{Introduction}

The smooth function extension problem is a classical problem and has been studied extensively in the literature from various view points. 
Some of the well-known results include the Uryshons Lemma, the Tietze Extension, and the Whitneys Extension Theorem (see, e.g, 
\cite{S, S1, BS, DF1, DF2, FDG, F5, Mc, Whitney}).

Inpainting was first introduced in \cite{BSBC} and then have been studied extensively by several authors (see, e.g., \cite{BBS}, \cite{CS}). 
Although smooth image inpainting is a smooth function extension problem, 
the most common approach in image inpainting so far is to use the solution to some PDE which are obtained from minimum-energy models as the recovered image. The most commonly used density function for these energy models is total derivation. 

In \cite{CS} by considering smooth inpainting as a smooth surface extension problem, the author studied methods for linear inpainting and cubic inpainting. 
Error bounds for these inpainting methods are derived in \cite{CS}. 
In \cite{C}, several surface completion methods have been studied. 
An optimal bound for the errors of the cubic inpainting method in \cite{CS} is given. 
Applications to smooth inpainting have also been discussed in \cite{C}. There, error bounds 
of completion methods are derived in terms of the radius of the domains on which  the functions are completed. In one of the methods in \cite{C}, the author proposed to 
use $p$-harmonic functions for smooth surface completion and smooth surface inpainting. Later in \cite{CM}, $p$-harmonic functions are also studied for smooth surface completion and inpaiting. The differences of the method using $p$-harmonic functions in \cite{CM} and \cite{C} are: the method in \cite{C} uses $\triangle^{i} u|_S$, $i=0,1,..,p-1$ as boundary data while the method in \cite{CM} uses $\frac{\partial^i}{\partial N^i}u|_S$, $i=0,1,...,p-1$ as boundary data to solve for a $p$-harmonic function. Here $u$ is the function to be inpainted and 
$S$ is the boundary of the inpainted region.

The goals of this paper are to implement and compare the performance of the two surface completion schemes in \cite{C} and \cite{CM}.  
In particular, we focus our study on smooth surface completion and smooth surface inpainting by biharmonic functions.  

\section{Surface completion by biharmonic functions}

Let $D$ be a simply connected region in $\mathbb{R}^2$ with 
$C^1$-boundary $S=\partial D$ and $d$ be the diameter of $D$. 
Let $u_0$ be a smooth function on any region containing $D$. Assume that 
$u_0$ is known on a neighbourhood outside $D$. The surface completion problem consists of finding a function $u$ on a region containing $D$ such that 
\begin{equation}
\label{eq2.1}
u = u_0\quad \text{outside}\quad D.
\end{equation}
There are several ways to construct the function $u$ 
so that \eqref{eq2.1} holds. For smooth surface completion, one often is interested in finding a sufficiently smooth function $u$ satisfying \eqref{eq2.1}. 

An application of smooth surface completion is in smooth image  inpainting. 
In smooth image inpainting, one has a smooth image $u_0$ which is known in a neighbourhood outside of a region $D$ while the data inside $D$ is missing. 
The goal of image inpainting is to extend the function $u$ over the region $D$ in such a way that the extension over the missing region is not noticeable 
with human eyes.  

In image inpainting, an inpainting scheme is said to be of linear order, or simply linear inpainting, if for
any smooth test image $u_0$, as the diameter $d$ of the inpainting region $D$ shrinks to 0, one has
$$
\|u - u_0\|_D = O(d^2)
$$
where $u$ is the image obtained from the inpainting scheme. 
Here $\|\cdot\|_D$ denotes the $L^\infty(D)$ norm. Here and throughout $f=O(g)$ if 
$\frac{|f|}{|g|}$ is bounded uniformly by some constant $C>0$. 

Note that harmonic inpainting, i.e., the extension found from the equation
$$
\triangle u = 0\quad \text{in} \quad D,\qquad u|_S = u_0|_S, 
$$
is a linear inpainting scheme (\cite{CS}). 

In \cite{CS} the following result for cubic inpainting is proved 
\begin{theorem}[Cubic Inpainting, Theorem 6.5 \cite{CS}] Let  $u_1$ 
be the harmonic inpainting of $u_0$. Let $u_\ell$ be a linear inpainting of $\triangle u_0$ on $D$ (not necessarily
being harmonic), and let $u(x)$ be defined by
\begin{equation}
\label{eq2.2}
u(x) = u_1(x) + u_2(x),\qquad x\in D,
\end{equation}
where $u_2$ solves the Poisson's equation
\begin{equation}
\label{eq2.3}
\triangle u_2 = u_\ell,\quad \text{in}\quad D, \quad u_2|_S = 0. 
\end{equation}
Then $u$
defines a cubic inpainting of $u_0$, i.e.,
$$
\|u - u_0\|_D = O(d^4).
$$
\end{theorem}

\begin{remark}{\rm 
If $u_\ell$ is the harmonic inpainting of $\triangle u_0$ in $D$, i.e., 
$u_\ell$ solves the equation
\begin{equation}
\label{eq2.4}
\triangle u_\ell = 0\quad \text{in} \quad D,\qquad u_\ell|_S = \triangle u_0|_S.
\end{equation}
Then the element $u$ defined by \eqref{eq2.2} is a biharmonic function which solves the following problem
\begin{equation}
\label{eq2.5}
\triangle^2 u = 0,\qquad \triangle u|_S = \triangle u_0|_S,\qquad 
u|_S = u_0|_S. 
\end{equation}

}
\end{remark}

In \cite{C}, this result is generalized to a multi-resolution approximation extension scheme 
in which the Laplacian is replaced by more general lagged diffusivity anisotropic operators. It is proved in \cite{C} that 
if $u$ solves the equation 
\begin{equation}
\label{eq2.6}
\triangle^n u = 0\quad \text{in}\quad D,\qquad 
\triangle^i u|_S = \triangle^i u_0|_S
, \quad i=0,1,...,n-1,
\end{equation}
then 
\begin{equation}
\label{eq2.7}
\| u - u_0 \|_D = O(d^{2n}).
\end{equation} 
A sharper error bound than \eqref{eq2.7} is obtained in \cite{C}. 

Equation \eqref{eq2.6} can be written as a system of Poisson's equations as follows
\begin{equation}
\label{eq2.9}
\left \{
\begin{split}
v_0 &= 0,\quad \text{in} \quad D,\\
\triangle v_i &= v_{i-1}\quad \text{in}\quad D, 
\quad v_i|_S = \triangle^{n-i} u_0|_S,\quad i=1,2,..,n,\\
u &= v_n.
\end{split}
\right.
\end{equation} 
Thus, the problem of solving \eqref{eq2.6} is reduced to the problem of solving Poisson's equations of the form
\begin{equation}
\label{eqd1}
\triangle u = f\quad \text{in}\quad D,\qquad u|_S = g. 
\end{equation}
Numerical methods for solving equation \eqref{eqd1} have been extensively studied in the literature.  

Note that the normal derivatives $\frac{\partial^i}{\partial N^i}u_0|_S$ 
are not presented in equation \eqref{eq2.6}. Thus, the extension using 
\eqref{eq2.6} may be not differentiable across the boundary $S$. 
To improve the smoothness of the extension across $S$, it is suggested to find $u$ from the equation
\begin{equation}
\label{eq2.8}
\triangle^n u = 0\quad \text{in}\quad D,\qquad \frac{\partial^i}{\partial N^i}u = \frac{\partial^i}{\partial N^i}u_0, \quad i=0,1,...,n-1. 
\end{equation}
It is proved in \cite{CM} that if $u$ is the solution to equation \eqref{eq2.8}, then \eqref{eq2.7} also holds.

Equation \eqref{eq2.8} cannot be reduced to a system of Poisson's equations as \eqref{eq2.6}. 
In fact, to solve \eqref{eq2.8} one often uses a finite difference approach which consists of finding discrete approximations to 
the operators $\triangle^n$ and $\frac{\partial^i}{\partial N^i}$, $i=1,...,n-1$. For `large' $n$, it is quite complicated, although possible, to obtain these approximations. 

As we can see from the above discussions, equation \eqref{eq2.6} 
is easier to solve numerically than equation \eqref{eq2.8}. 
However, scheme \eqref{eq2.6} does not use any information about the normal derivatives of the surface on $S$. 
Thus, the extension surface, obtained from scheme \eqref{eq2.6}, may not be smooth cross the boundary $S$. On the contrary, equation \eqref{eq2.8} uses normal derivatives as boundary data and, therefore, is expected to yield better results than scheme \eqref{eq2.6} does.  

In the next section we will implement and compare the two surface completion schemes using equations \eqref{eq2.6} and \eqref{eq2.8}. 
In particular, we focus our study on biharmonic functions which are solutions to \eqref{eq2.6} and \eqref{eq2.8} for $n=2$.

\section{Implementation}

Let us discuss a numerical method for solving the equation
\begin{equation}
\triangle^2 u = 0\quad\text{in}\quad D,\quad \triangle u|_S = f,\quad u|_S = g. 
\end{equation}
To solve this equation, one often defines $v=\triangle u$ and solves 
for $u$ from the following system
\begin{equation}
\label{eq2}
\triangle v = 0,\quad v|_S = f,\qquad \triangle u = v,\quad u|_S = g. 
\end{equation}
Thus, the problem of solving \eqref{eq2} is reduced to the problem of solving the following Poisson's equation
\begin{equation}
\label{eq3}
\triangle u = f_1\quad \text{in}\quad D,\qquad u|_S = g_1. 
\end{equation}
To solve equation \eqref{eq3}, we use a 5-points finite difference scheme to approximate the Laplacian operator. This 5-points scheme is based on the following well-known formula
\begin{equation}
\label{eq4}
\triangle_5 u := \frac{1}{h^2}
\begin{pmatrix}
&1 &\\
1& - 4& 1\\
&1&\\
\end{pmatrix}u =\triangle u + O(h^2).
\end{equation}
Here, $h$ is the discretization step-size. 
This scheme is well-defined at points $P$, whose nearest neighbours are interior points of $D$. If $P$ has a neighbour $Q\in \partial D$, then we use a stencil of the form
\begin{equation}
\label{eq4i}
\triangle_5 u(P) := \frac{1}{h^2}
\begin{pmatrix}
1 &\\
 - 4& 1\\
1&\\
\end{pmatrix}u(P) + \frac{u(Q)}{h^2}=\triangle u(P)+O(h^2).
\end{equation}
In the above formula $Q$ is the nearest neighbour to the left of $P$. Similar formulae 
when $Q$ is the nearest neighbour on the top, on the bottom, and to the right of $P$ can be obtained easily.

In our experiments, we choose 
$D$ as a square  
and the solution $u$ on the computation grid is presented as a vector.   
Using the 5-points finite difference scheme above equation \eqref{eq3} is reduced to the following algebraic system
\begin{equation}
\label{eq6}
Au = (f_1 - \frac{\tilde{g}_1}{h^2}),
\end{equation}
where $A$ is the 5-points finite difference approximation to the Laplacian and $\tilde{g}_1$ is a vector containing boundary values of $u$ on $S$ at suitable entries. 
The matrix $A$ is is a tridiagonal matrix, that is, all non-zero elements of $A$ lie on the main diagonal, and the first diagonals above and below the main diagonal.  
The matrix $-h^2A$ can be obtained by the function $delsq$ available in Matlab.

Let us discuss a numerical method for solving the equation  
\begin{equation}
\label{eq7}
\triangle^2 u = 0,\quad u_N|_S = f,\quad u|_S = g. 
\end{equation}
For a discrete approximation to the bilaplacian we use a 13-points finite difference scheme  which is based on the following formula (see \cite{B})
\begin{align}
\triangle_{13}^2u = 
\frac{1}{h^4}
\begin{pmatrix}
&&1&&\\
&2&-8&2&\\
1&-8& 20&-8&1\\
&2&-8&2&\\
&&1&&
\end{pmatrix}
u = \triangle^2 u + O(h^2).
\end{align}
This stencil is well-defined for a grid point $P$ if all its nearest neighbours are in 
the interior of the domain. 
If $P$ has a neighbour $Q\in\partial D$ and $Q$ is on the left of $P$, then we use the following formula 
\begin{align}
\triangle_{13}^2u(P) :&= \frac{1}{h^4}
\begin{pmatrix}
&1&&\\
2&-8&2&\\
-8& 20&-8&1\\
2&-8&2&\\
&1&&
\end{pmatrix}
u(P)+\frac{u_N(Q)}{h^3} + \frac{u(Q)}{h^4} \nonumber\\
&= \triangle^2 u + O(h^{-1}).
\end{align}

Using the above finite difference scheme equation \eqref{eq7} is reduced to a linear algebraic system of the form $Au = b$ where $A$ is a five-diagonal matrix. Numerical solutions to $u$ on the grid can be obtained by solving this linear algebraic system. 
\medskip

Before we proceed with numerical experiments we need:

\subsection{Quantative comparisons}

It is constructive to provide quantitative correlations between original and processed images and in particular code to compare figures such as those below. 
In order to calculate these required correlations (and many are provided), we refer the reader to a free access code for our method  in a sckit-image processing in python unit competely at the disposal of the reader which readily provides  quantative correlations---in particular for the figures below. 
The  sckit-image-image processing package is at:
\medskip

 http://scikit-image.org/ and is a collection of open access algorithms for image processing with peer-reviewed code. 
\medskip

For our method, see:
\medskip

http://scikit-image.org/docs/stable/api/skimage.restoration.html\newline
?highlight=biharmonicskimage.restoration.inpaint-biharmonic .
\medskip

Moreover, for the benefit of the reader, many comparison methods with associate code are given in this unit,  see  below for some, but see full package for a longer list with references.

\begin{itemize}
\item Denoise-bilateral, see \cite{TM}
\item Denoise-nl-means, see \cite{BCM,DCCOJ,F}
\item Denoise-tv-bregman, denoise-tv-chambolle, denoise-wavelet, estimate-sigma, inpaint-biharmonic (the present paper), nl-means-denoising, richardson-lucy, unsupervised-wiener, unwrap-phase, Wiener-Hunt deconvolution.
\end{itemize}

\section{Numerical experiments}

\subsection{Smooth surface completion}

Let us first do some numerical experiments with smooth surface completion. 
In our experiments, we compare numerical solutions from the three surface completion methods: the method by harmonic functions, the method by biharmonic functions in \cite{CS} and \cite{C}, and the method with biharmonic functions in \cite{CM}.

The function $u$ to be completed in our first experiment is
$$
u(x,y) = xy + x^2(y+1),\qquad x,y \in D=[-1,1].
$$
Note that this function is a biharmonic function. 
The domain $D$ is discretized by a grid of size $(n+1)\times (n+1)$ points. 
In the first experiment  we used $n=50$. 

In our numerical experiments, we denote by $u_H$ the extension by a harmonic function, by $u_L$ the biharmonic extension from \cite{C}, and by $u_N$ the biharmonic extension from \cite{CM}.   

\begin{figure}[!h!tb]
\centerline{%
\includegraphics[scale=0.93]{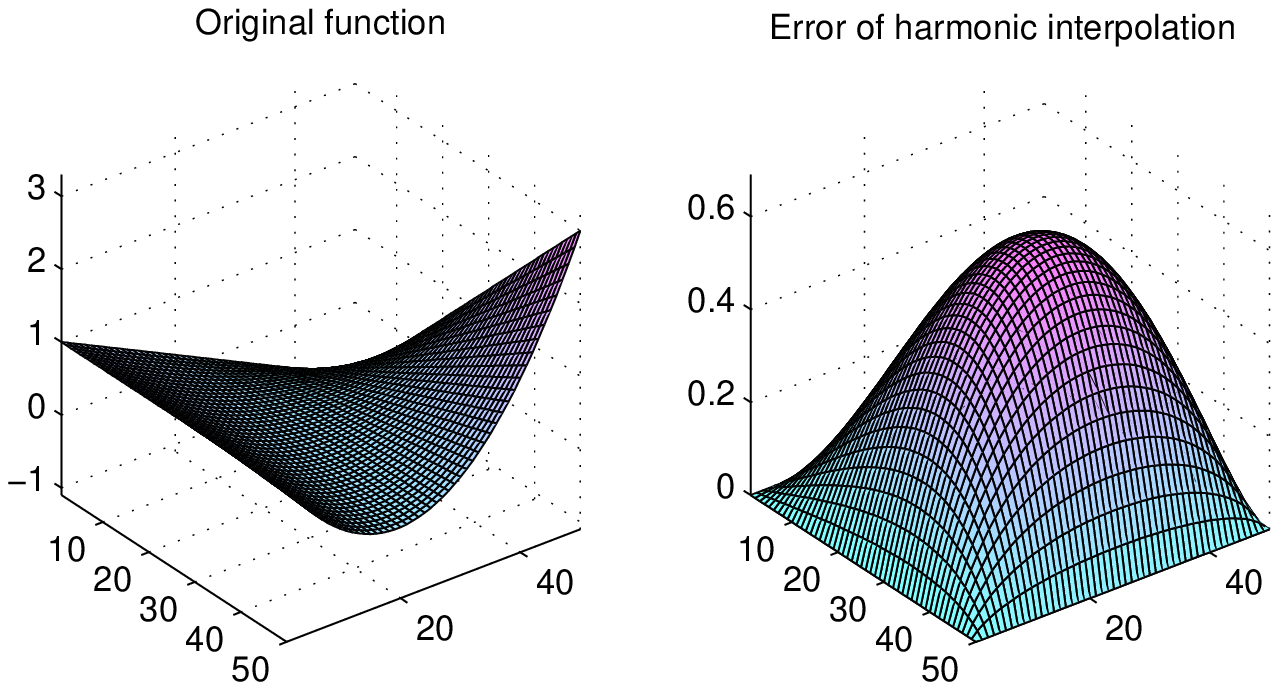}}
\caption{}
\label{example1a}
\end{figure}

\begin{figure}[!h!tb]
\centerline{%
\includegraphics[scale=0.93]{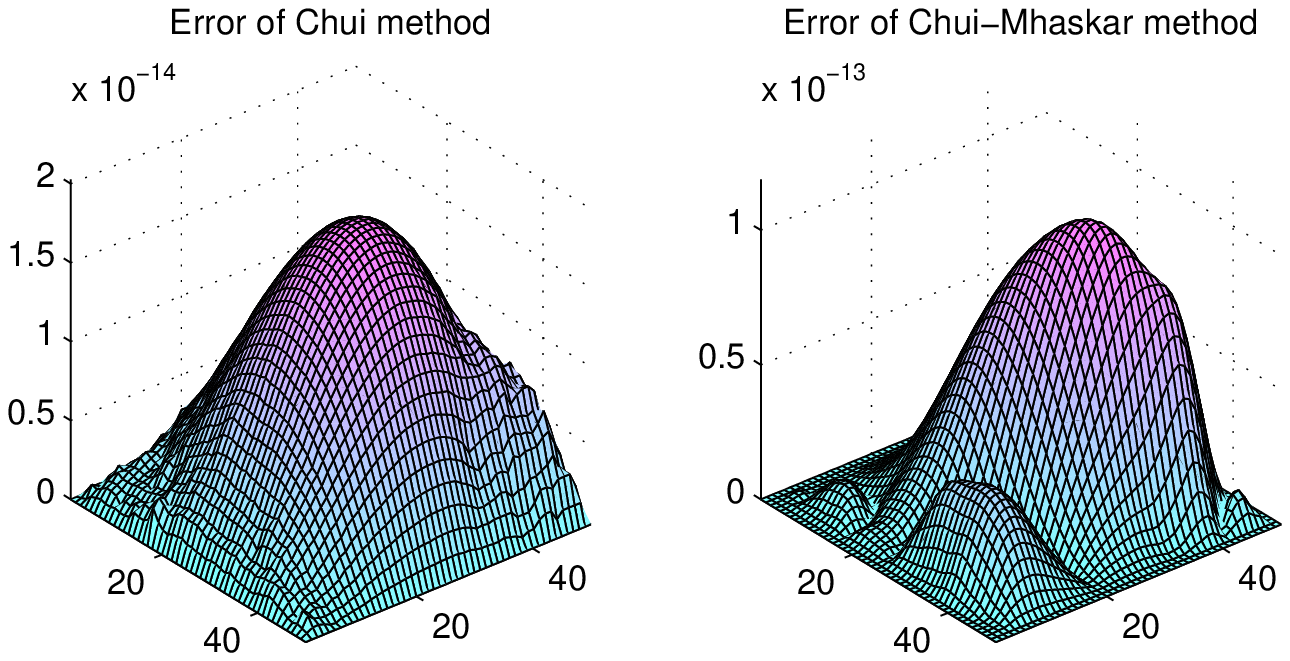}}
\caption{}
\label{example1b}
\end{figure}

Figure \ref{example1a} plots the original function and the error of 
reconstruction by a harmonic function. Figure \ref{example1b} plots the errors of 
surface reconstructions by biharmonic functions from \cite{C} and \cite{CM}.

From Figures \ref{example1a} and \ref{example1b}, one can see that 
the biharmonic reconstructions from \cite{C} and \cite{CM} are much better than the 
reconstruction by harmonic functions. 
The method in \cite{C} in this experiment yields numerical results with accuracy a bit 
higher than the method in \cite{CM}. However, this doesn't imply that the method in \cite{C} is better in term of accuracy than the method in \cite{CM}. The condition number of $A$, the finite difference approximation to the bilaplacian in this experiment is larger than that of the finite difference approximation to the Laplacian. Due to these condition numbers, the algorithm using 
the method in \cite{C} yields results with higher accuracy then the algorithm using 
the method in \cite{CM}. As we can see in later experiments, the method in \cite{CM} often gives better results than the method in \cite{C}. 
The conclusion from this example is that both the methods in \cite{C} and \cite{CM} 
yield numerical solutions at very high accuracy. The harmonic reconstruction in this experiment is not very good. This comes from the fact that the function to be reconstructed is not harmonic. 

In the next experiment the function to be reconstructed is chosen by
\begin{equation}
\label{eq3.9}
u(x,y) = \frac{(1+\cos(x))(1+\cos(y))}{4},\qquad x,y\in[-1,1].
\end{equation}
This function $u(x,y)$ is not a biharmonic function.

\begin{figure}[!h!tb]
\centerline{%
\includegraphics[scale=0.75]{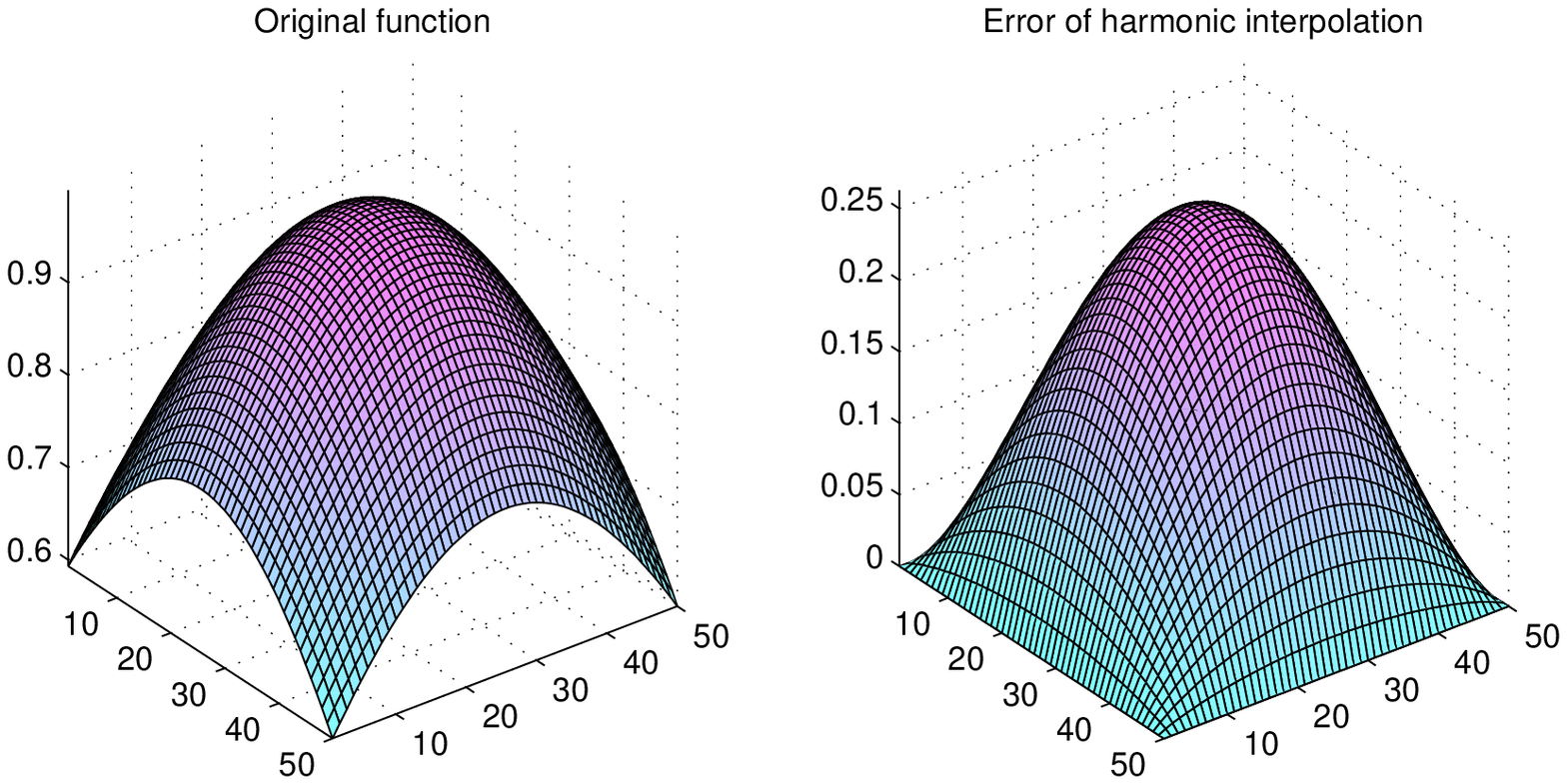}}
\caption{}
\label{example2a}
\end{figure}

\begin{figure}[!h!tb]
\centerline{%
\includegraphics[scale=0.75]{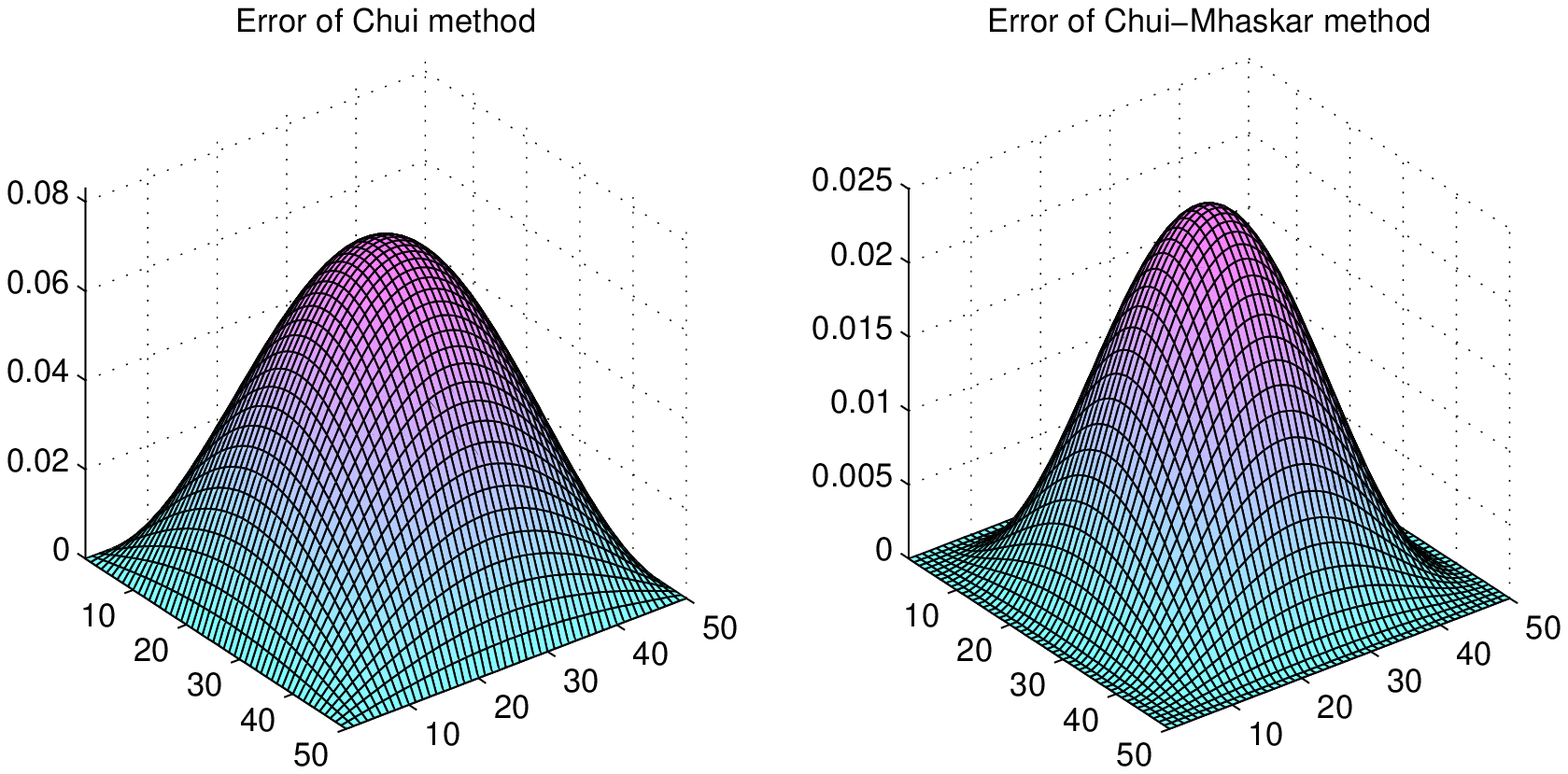}}
\caption{}
\label{example2b}
\end{figure}

Figure \ref{example2a} and \ref{example2b} plot the errors of harmonic and biharmonic reconstructions. From these figures it is clear that the method \cite{CM} yields the best approximation. The harmonic reconstruction is the worse amongst the 3 methods in this experiment. 

\subsection{Image inpainting}

Let us do some numerical experiments with image inpainting. 

\begin{figure}[!h!tb]
\centerline{%
\includegraphics[scale=0.93]{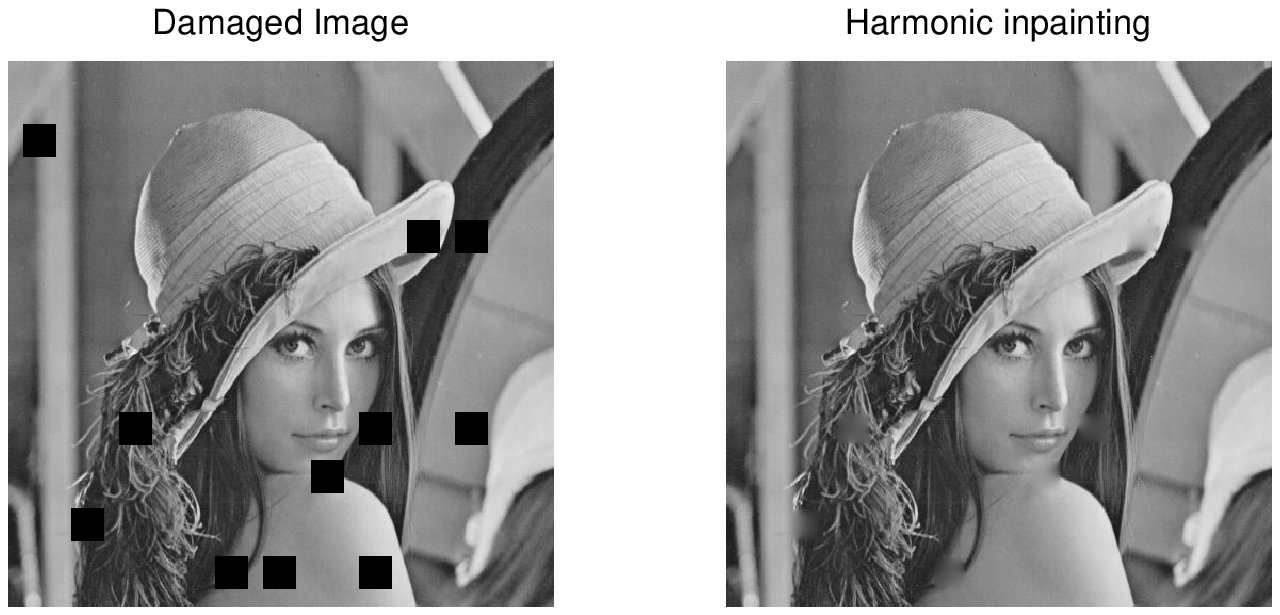}}
\caption{Damaged image and restored image by harmonic functions.}
\label{lena1}
\end{figure}

\begin{figure}[!h!tb]
\centerline{%
\includegraphics[scale=0.93]{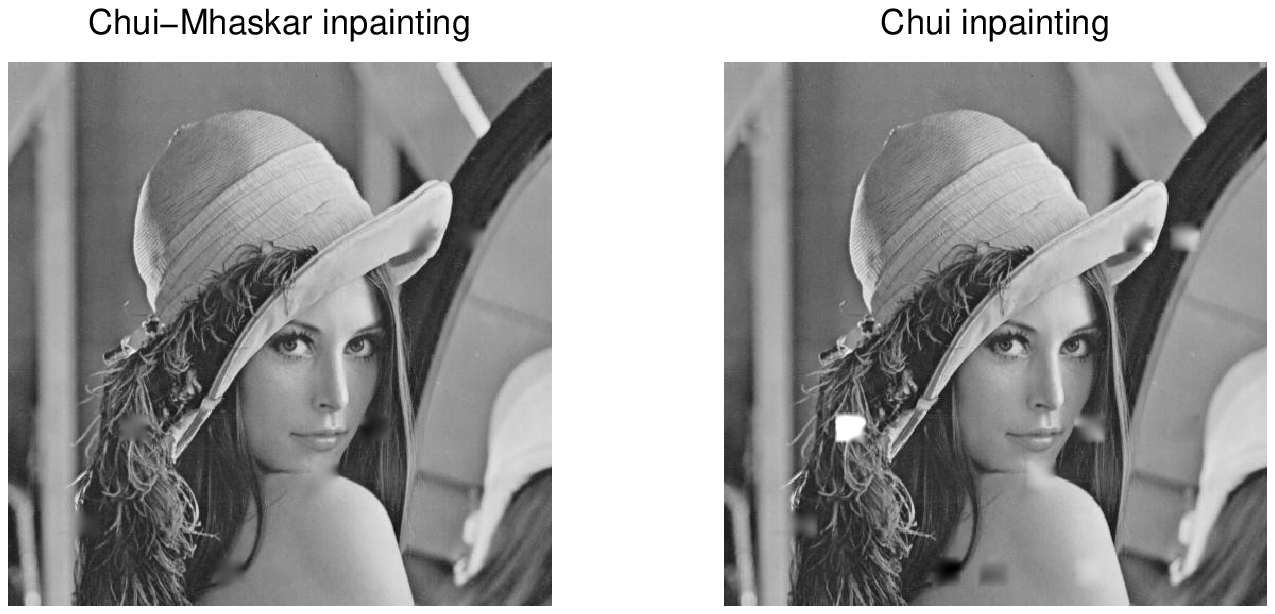}}
\caption{Restored image by biharmonic functions.}
\label{lena2}
\end{figure}

Figure \ref{lena1} plots a damaged image and a reconstructed image by harmonic functions. Figure \ref{lena2} plots restored images by biharmonic functions following the methods from \cite{C} and \cite{CM}. 

It can be seen from Figures \ref{lena1} and \ref{lena2}, that the biharmonic extension method from \cite{CM} yields the best reconstruction. Although the biharmonic extension method from \cite{C} is better than harmonic extension in our experiments with smooth surface completion, it is not as good as the harmonic extension in this experiment. This is understandable since our image contains edges and is not a smooth function. It can be seen from the restored image by the method in \cite{C} in Figure \ref{lena2} that the reconstruction may not be smooth, or even not differentiable, across the boundary.

\medskip
\medskip

S.B. Damelin: Mathematical Reviews, The American Mathematical Society, 416 Fourth Street, Ann Arbor, MI 48104. email: damelin@umich.edu 
\medskip

N. S. Hoang: Department of Mathematics, University of Oklahoma,
Norman, OK 73019-3103, USA: email: nhoang@math.ou.edu.

\medskip

The authors declare that there is no conflict of interest regarding the publication of this paper.
\medskip

\section{Appendix: figures and tables.}

Table \ref{table1} presents numerical results for the function $u$ defined by \eqref{eq3.9} and on the domain $D_i=[-2^{-i},2^{-i}]\times [-2^{-i},2^{-i}]$, $i=\overline{0,9}$. The diameter of $D_i$ is $d_i=2^{1-i}$. From Table \ref{table1}, one can see that the harmonic reconstruction has an order of accuracy of $2$ while the biharmonic reconstruction methods have an order of accuracy of $4$. This agrees with the theoretical estimates in \cite{CS}, \cite{C}, and \cite{CM}. 

\begin{table}[ht]
\caption{Results for $D_i=[-2^{-i},2^{-i}]\times [-2^{-i},2^{-i}]$ for $i=\overline{0,9}$}
\label{table1}
\centering
\small
\begin{tabular}{@{  }c@{\hspace{3mm}}|c@{\hspace{2mm}}|c@{\hspace{2mm}}|c@{\hspace{2mm}}|}
\hline
i & $\log_2\|u_H-u\|_\infty$& $\log_2\|u_L-u\|_\infty$ & $\log_2\|u_N-u\|_\infty$\\
\hline
    0&    2.36&    0.37&   -1.46\\
    1&    0.50&   -3.48&   -5.34\\
    2&   -1.46&   -7.44&   -9.32\\
    3&   -3.45&  -11.44&  -13.31\\
    4&   -5.45&  -15.43&  -17.31\\
    5&   -7.45&  -19.43&  -21.31\\
    6&   -9.45&  -23.43&  -25.31\\
    7&  -11.45&  -27.43&  -29.31\\
    8&  -13.45&  -31.43&  -33.34\\
    9&  -15.45&  -35.46&  -38.54\\
\hline
\end{tabular}
\end{table}

\begin{figure}[!h!tb]
\centerline{%
\includegraphics[scale=0.95]{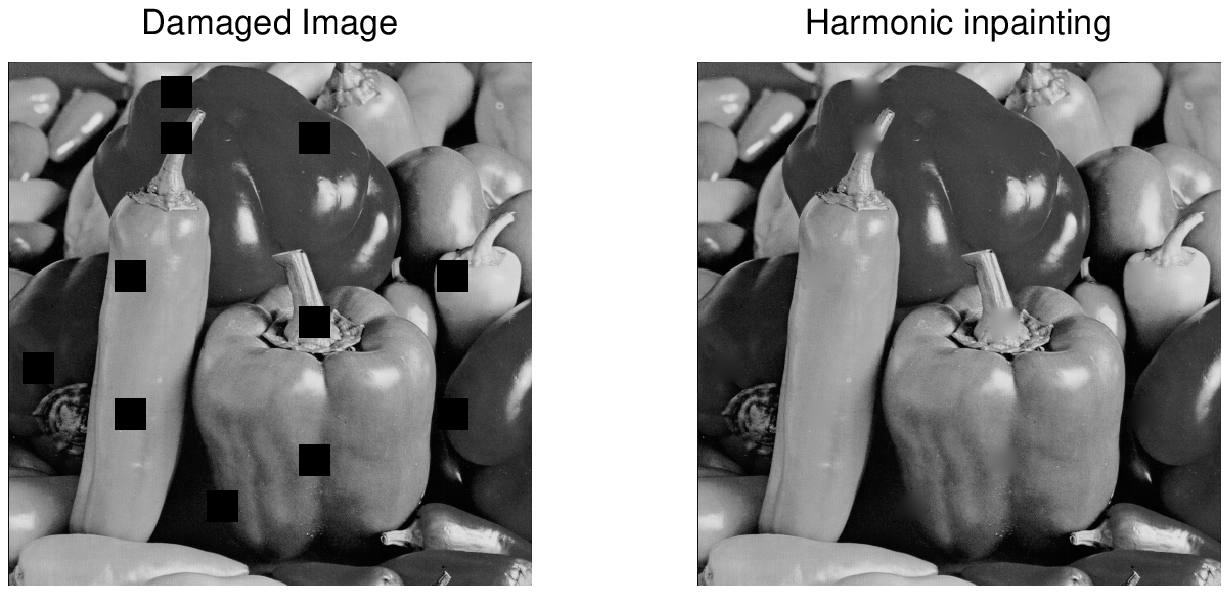}}
\caption{Damaged image and restored image by harmonic functions.}
\label{peppers1}
\end{figure}

\begin{figure}[!h!tb]
\centerline{%
\includegraphics[scale=0.95]{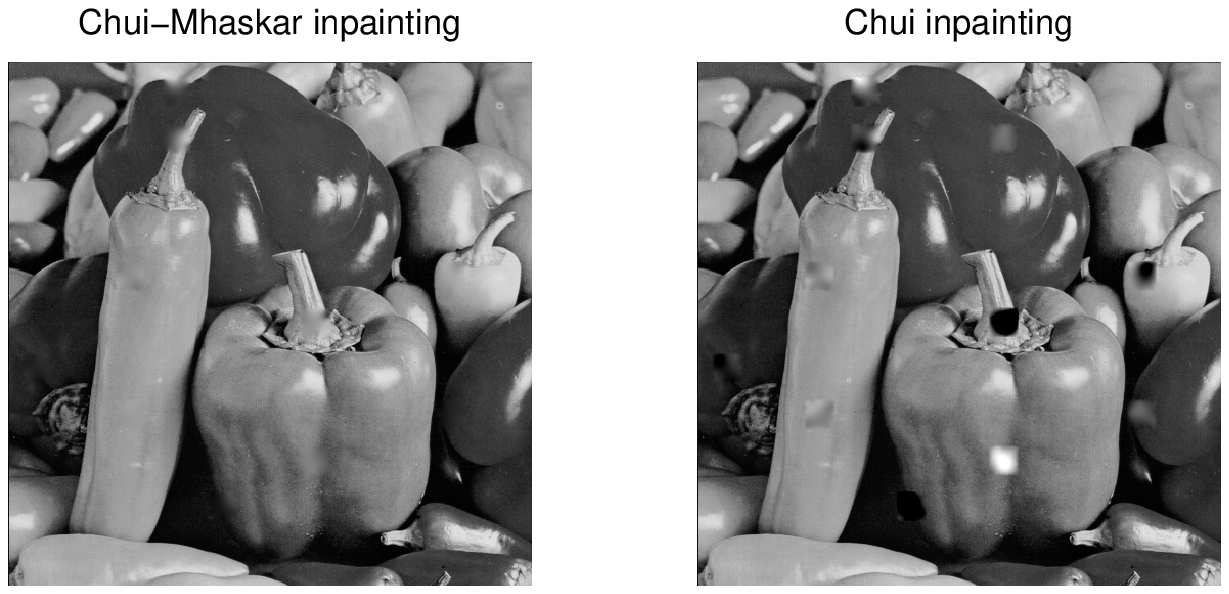}}
\caption{Restored image by biharmonic functions.}
\label{peppers2}
\end{figure}

Figures \ref{peppers1} and \ref{peppers2} 
plot a damaged picture of peppers and reconstructed images by the 3 methods. 
From these figures one gets the same conclusion as in the previous experiment. 
The biharmonic reconstruction in \cite{CM} yields the best restoration while the biharmonic reconstruction in \cite{C} yields the worst reconstruction. 

\end{document}